\documentstyle[12pt,amstex]{article}
\title{Hopf-Galois extensions with central invariants.}
\author{Dmitriy Rumynin\footnote{University of Massachusetts at Amherst}}
\date{September 22, 1997}

\newtheorem{theorem}{Theorem}
\newtheorem{prop}[theorem]{Proposition}
\newtheorem{lemma}[theorem]{Lemma}
\newtheorem{dfn}[theorem]{Definition}

\newtheorem{ldf}[theorem]{Theorem-Definition}

\begin{document}
\maketitle

\begin{abstract}
We study a class of algebra extensions which usually
appear in the study of restricted Lie algebras or various
quantum objects at roots of unity.
\end{abstract}

The present paper was inspired by the theory of nonrestricted 
representations
of restricted Lie algebras 
and the theory of quantum groups at roots of unity
where  algebras are usually finitely generated modules
over their centers.
Our objective is to demonstrate that all these theories
admit a unifying approach.
We use the concept of a Hopf-Galois extension
with central invariants
to treat these phenomena from a general point of view.
We discuss extensions of algebras
 in Section 1. We define Hopf-Galois extensions
with central invariants early in Section 2. Then
we study them locally: i.e. we consider localizations
at points of the prime spectrum of the subalgebra of invariants.
The localizations form a vector bundle of not necessarily
commutative algebras on a scheme.
The fibers of the bundle are
finite-dimensional Frobenius algebras and their irreducible
representations coincide with irreducible representations
of the Hopf-Galois extension we study.
We describe some known examples in which the developed
formalism takes place in Section 3.
It seems that further attempts to study these extensions
should be based on exploring the inherent geometry of the 
situation.
The author is grateful to J. Humphreys for fruitful
discussions and the proofreading of the manuscript.

\section{Preliminary facts on algebra extensions.}
We are working over a ground field ${\bf k}$ from
now on. 
Concerned with an algebra $U$ (associative with unity),
 we consider  a finitely-generated
subalgebra $O$ of $Z$, the center of $U$, 
such that the module $U_O$ is finitely-generated.
Such extensions were first studied in \cite{cur}. 
In this section we collect various technical results needed later.

Under the assumptions we make $U$ is a PI-algebra. One should use
\cite{pro} as a general reference. The following proposition easily follows
from the Hilbert basis theorem.

\begin{prop}
  $Z_O$ is a finitely generated module and $Z$ is a finitely generated
  algebra.
\end{prop}

The following is a version of Hilbert's Nullstellensatz for $U$.
\begin{theorem}
  A simple $U$-module is finite dimensional.  Furthermore, if $\rho$
  is a finite dimensional representation then $\rho (Z)$ and $\rho
  (O)$ are finite extensions of $\bf k$.
\end{theorem}

The next proposition is similar to \cite[Theorem 2.7.4]{pro}.

\begin{prop} Simple $U$-modules are isomorphic if and only if the kernels
  of the corresponding representations are equal.
\end{prop}

Proposition 3 shows that assigning the kernel to an irreducible
representation is a bijection from the set of simple $U$-modules ${
  \frak Rep }U$ onto the set of primitive ideals. The latter has the
structure of a topological space with respect to the Jacobson
topology. The central characters $\Omega$ yield a sequence of maps:
\begin{equation}
{ \frak Rep }U \cong \mbox{Prim }U
\stackrel{\Omega_Z}{\longrightarrow} \mbox{Max }Z
\stackrel{\Omega_{Z/O}}{\longrightarrow} \mbox{Max }O \label{1}
\end{equation}
 We
denote the composition $\Omega_{Z/O} \circ \Omega_{Z}$ by $\Omega_O$.
The following lemma is a non-commutative variation of a well-known
result in algebraic geometry \cite{cur}.
 
\begin{lemma}
  $\Omega_Z$ and $\Omega_{Z/O}$ are continuous onto maps.
\end{lemma}

The finitely generated module $U_O$ defines a coherent sheaf of
algebras on Spec$\: O$. Denoting by $O_{(\eta)}$ the local ring of
algebraic functions at $\eta \in \mbox{Spec} \: O$, $m_{\eta}$ its
maximal ideal, and $K_\eta$ the quotient field $O_{(\eta)} / m_\eta$,
we introduce three algebras:
\begin{equation}U_{(\eta)}=U \otimes_{O} O_{(\eta)} \mbox{ is an } O_{(\eta)} \mbox{-algebra;}
\label{2} \end{equation}
$$U_\eta =
U_{(\eta)} \otimes_{O_{(\eta)}} K_{\eta} \mbox{ is a }K_\eta \mbox{-algebra;}$$
$$U_{[\eta]}=U / U \eta \mbox{ is an } O/ \eta \mbox{-algebra.}$$

If the module $U_O$ is projective
the sheaf becomes an algebraic vector bundle. 
The following theorem holds \cite[Theorem 4.27]{mcd}:
\begin{theorem}
If $U_O$ is projective then the dimension function
$\chi \mapsto \mbox{\em dim}_{K_{\chi}}
 \: U_\chi$ is a continuous map from
{\em Spec }$O$ to $\Bbb Z$.
\end{theorem}
We would like to understand 
 $O$-bilinear forms $s:U \times U \longrightarrow O$.
Let a $K_\eta$-bilinear form $s_{\eta}:U_{\eta} \times U_{\eta} 
\longrightarrow K_{\eta}$ be the specialization of $s$  to a point $\eta \in
\mbox{Spec} \: O$.
The proof of the following theorem has the same main idea
as in the theorem about continuity of the dimension function.

\begin{theorem} Let the module $\,_OU$ be projective. 
The set $\{ \eta \in \mbox{\em Spec} \: O \mid
s_{\eta} \mbox{ is non-degenerate} \}$  is open.
\end{theorem}

{\bf Proof.}
$s$ defines an $O$-module map
$\hat{s}: \,_OU \longrightarrow (\,_OU)^\ast$.
Then $s_{\eta}$ is non-degenerate outside the union of supports
of the kernel and the cokernel $\hat{s}$. We will give a precise
argument which works for not necessarily finitely generated $O$.

The specialization $s_{\eta}$ is non-degenerate if and only if
$\hat{s}_{\eta}: \,_{K_{\eta}}U_{\eta} \longrightarrow 
(\,_{K_{\eta}}U_{\eta})^\ast$ is 
\linebreak
an isomorphism.
It is equivalent the localization 
$\hat{s}_{(\eta )}
: \,_{O_{(\eta )}}U_{(\eta )} \longrightarrow 
(\,_{O_{(\eta )}}U_{(\eta )})^\ast$
being an isomorphism. Indeed,
the projective module $U_{(\eta )}$ over the local ring $O_{(\eta )}$
must be free which allows us to use the determinant to check whether
a map is an isomorphism. Clearly, det$\, \hat{s}_{(\eta)} +m_{\eta}=
\mbox{det} \, \hat{s}_{\eta}$. Thus, det$\, \hat{s}_{(\eta)}$ is invertible
in $O_{(\eta )}$ if and only if 
det$\, \hat{s}_{\eta }$ is invertible in $K_{\eta}$.

The map $\hat{s}$ is a part of the exact sequence
$0 \rightarrow \:_OA \rightarrow \,_OU \stackrel{\hat{s}}{\longrightarrow} 
(\,_OU)^\ast \rightarrow \,_OB \rightarrow 0$.
If $s_\eta$ is degenerate for each $\eta$ then we have nothing to prove:
the empty set is open. Let us pick $\eta$ with non-degenerate $s_\eta$.
Now $\,_OB$ is finitely generated as a quotient of a finitely generated module.
 Let $b_1, \ldots , b_n$ be the generators.
Since $\,_OB_{(\eta )}=0$ we can find elements $x_i \in O$ such that
$x_i b_i =0$ and $x_i \notin \eta$. Let $x=x_1 \cdots x_n$. We have obtained 
the exact
sequence of localizations:
$0 \rightarrow \:_OAx^{-1} \rightarrow \,_OUx^{-1} 
\stackrel{\hat{s}x^{-1}}{\longrightarrow}
(\,_OU)^\ast x^{-1} \rightarrow \,_OBx^{-1} \cong 0$.
The modules $\,_OUx^{-1}$ and $(\,_OU)^\ast x^{-1}$ are projective.
Thus, the latter sequence is split and $\,_OA x^{-1}$ is finitely
generated. By the argument as above we can find $y \in O \setminus \eta$
such that $(Ax^{-1})y^{-1} =0$. This implies that $(\hat{s}x^{-1}
)y^{-1}$ is isomorphism. Therefore, $s_\eta$ is non-degenerate
for $\eta$ from the open set $\{ I \mid x \notin I \mbox{ and }
y \notin I \}$. $\Box$

If $U_O$ is projective the theory of deformations \cite{ger} may be used 
to get some information about the bundle.
Let $X$ be the points of Spec$\; O$ over a field $F$.
\begin{prop} Let $\,_OU$ be projective. 
The set $\{ \chi \in X \mid U_{\chi} \mbox{ is
a separable } K_{\chi} \mbox{-algebra} \}$ is open in
the Zariski topology on $X$.
\end{prop}
{\bf Proof.} The algebra $U_\chi$ is separable if and only if
so is $\tilde{U}_{\chi}= U_{\chi} \otimes_{K_{\chi}} F$.
Let $\eta \in X$ be a point such that $\tilde{U}_{\eta}$ is separable.
We denote the dimension of $\tilde{U}_{\eta}$ over $F$ by $n$.
By Theorem 5 there exists an open neighborhood $W$ of $\eta$ such that
for each $\chi \in W$ the dimension of $\tilde{U}_{\chi}$ is $n$.

Let $\mbox{Ass}_n \subseteq (F^n)^{\ast} \otimes (F^n)^{\ast} \otimes F^n$ 
consist of all tensors defining a structure of an associative
algebra on $F^n$. Thus we get a map $\psi :W \longrightarrow
\mbox{Ass}_n/ \mbox{GL}_n(F)$. The latter may be either a topological space 
with $\psi$
being a continuous map or an algebraic variety with a rational map $\psi$.
It is proved in \cite{ger} that the set 
$V_A = \{ B \in \mbox{Ass}_n \mid A\otimes_F \overline{F}
\cong B \otimes_F \overline{F}$ is open 
for each separable algebra $A$. 
$V_A/\mbox{GL}_n(F)$ is open in $\mbox{Ass}_n/ \mbox{GL}_n(F)$
since $V_A$ is saturated. 
Therefore, the set 
$\psi^{-1} (V_{\tilde{U}_{\eta}}/ \mbox{GL}_n(F))$ is an open 
neighborhood of $\eta$ with separable algebras. $\Box$

We should notice that we have proved a little more. All algebras
on some open subset of $X$ are separable and isomorphic
over the algebraic closure.
In particular, if $X$ is irreducible and $F$ is algebraically
closed then all separable
algebras of the bundle must be isomorphic.
However, we have not proved that the open set is not empty.
For instance, if $H$ is not semisimple and $U$ is just
a tensor product $O \otimes H$ then the open set on which
$U_{\chi}$ is separable is empty.

Another important restriction one may impose is $U$ being prime,
i.e. a product of two proper ideals is proper.
It implies that both $Z$ and $O$ are domains. Moreover,
the set $O^{\circ}$ of non-zero elements of $O$
is a set of non-zerodivizors of $U$.
The localizations with respect to $O^{\circ}$ behave nicely:  $Q(Z)$ is a field
and $Q(U)$  is a central
 $Q(Z)$-algebra. In particular,
the dimension of $Q(U)$ over $Q(Z)$ is square. The following
theorem may be found in  \cite[Corollary 5.3.2]{pro}.
\begin{theorem} If $U$ is prime then
the map $\Omega_Z$ (see~\ref{1}) is a bijection on the complement of a closed set
of smaller dimension.
\end{theorem}

\section{Hopf-Galois extensions.}

We use the standard notation 
for the structure maps of a Hopf algebra $H$: $\Delta, \varepsilon , S$
denote the comultiplication, the counity, and the antipode.
We write $\Delta (x) = x_1 \otimes x_2$ keeping summation in mind.
All tensor products are over {\bf k} unless otherwise
indicated.

\subsection{Basic properties of Hopf-Galois extensions.}

Let $H$ be a finite dimensional Hopf algebra. $U$ is an $H$-comodule algebra
if there is a map
$\rho :U \longrightarrow U \otimes H \ , \ \ \ \rho (x)=x_0 \otimes x_1$
such that $x_0 \varepsilon (x_1) =x$, $x_0 \otimes \Delta (x_1) = \rho (x_0) 
\otimes x_1$, and $\rho (xy)= \rho (x) \rho (y)$ for all $x,y \in U$.
If $O= U^H = \{ x \in U \mbox{ such that } \rho (x)=x \otimes 1 \}$ is a subalgebra
of invariants one calls $U \supseteq O$ an $H$-extension.
An $H$-extension is an extension with central invariants
if $O$ is contained in the center of $U$. An $H$-extension is Hopf-Galois
(or specifically $H$-Galois) 
if the canonical map can:$U \otimes_O U \longrightarrow U \otimes H$
defined by can$(x \otimes_O y)=(x \otimes 1)(y_0 \otimes y_1)$ is onto
\cite{krt,mon}.

We study $U \supseteq O$, a Hopf-Galois extension with central invariants,
such that $O$ is a finitely generated algebra
throughout this section. Here we list some important properties of such extensions.

\begin{prop} \cite[Theorem 1.7]{krt}.
 $U_O$ is a projective finitely generated module.
\end{prop}

\begin{ldf} \cite[7.2.2;8.2.4]{mon}.
The extension is called cleft if  the following equivalent
conditions hold. 

1. $U$ is isomorphic to a crossed product $O \#_{\sigma} H$ as an algebra.

2. There exists a convolution invertible right $H$-comodule map
$\gamma : H \longrightarrow U$.

3. There is a linear map from $U$ to $O \otimes H$ which is an isomorphism
of left $O$-modules and  right $H$-comodules.
\end{ldf}

\begin{lemma} Let $A \subseteq B$ be a cleft $H$-Galois extension such that
$A$ is commutative. $A$ is a central subalgebra if and only if $H$ acts
trivially on $A$.
\end{lemma}

{\bf Proof.} Let $\gamma : H \longrightarrow U$ be a splitting map
(i.e. a convolution invertible right $H$-comodule map).
If $A$ is central for every $a \in A$, $h \in H$ we obtain  
$h \cdot a =  \gamma (h_1) a \gamma^{-1} (h_2) =
 a \gamma (h_1) \gamma^{-1} (h_2) = a \varepsilon (h)$
Conversely, if the action is trivial we may think of $B$ as a crossed
product. 
$$(a \# h)( b \# 1)= a (h_1 \cdot b) \sigma (h_2 \otimes 1) \# h_3 =
ab \# h = (b \# 1)( a \# h) \ \ \ \ \ \Box$$

We have intentionally not introduced the crossed products.
An anxious reader can find them in \cite{mon}.
The last lemma indicates that we are interested only in
crossed products with the trivial action, i.e. in twisted products.
In the next section
we mostly reformulate the well-known properties of crossed products 
for twisted products although Theorem 14 has no crossed
products counterpart.

\subsection{Twisted products.}

Let $R$ be a commutative ${\bf k}$-algebra and $\sigma : H \otimes H \longrightarrow R$
be a linear map. Let $R_{\sigma}[H]$ be a vector space
$R \otimes H$ with an algebra structure defined by the formula
\begin{equation} a \otimes g \cdot b \otimes h =  ab \sigma (h_1,g_1) \otimes h_2g_2
\ , \ \ \ a,b \in R , g,h \in H \label{3} \end{equation}
It may not be a structure of an associative algebra in general.
The following lemma is straightforward \cite[Lemma 7.12]{mon}.

\begin{lemma}
The formula~\ref{3} defines a structure of an associative algebra
with an identity element $1 \otimes 1$ if and only if
for each $g,h,t \in H$
$$ \sigma (h_1 \otimes g_1) \sigma (h_2 g_2 \otimes t)=  
\sigma (g_1 \otimes t_1) \sigma (h \otimes g_2t_2)$$
$$ \sigma (h \otimes 1) = \sigma (1 \otimes h) = \varepsilon (h) 1$$
\end{lemma}

If $\sigma$ is invertible with respect to the convolution \cite{mon}
and satisfies the conditions of Lemma 12 we call it a cocycle,
although a corresponding cochain complex has been constructed
only for cocommutative $H$.
We use the term ``twisted product'' for $R_{\sigma}[H]$ 
if $\sigma$ is a cocycle.
 The following proposition
is a reason to introduce cohomological equivalence of cocycles \cite[7.3.4]{mon}.

\begin{prop}
Two cocycles $\sigma$ and $\tau$ produce isomorphic twisted products
if and only if there exists a linear convolution-invertible map
$u : H \longrightarrow R$ such that
$\tau (h \otimes g) = u^{-1} (g_1) u^{-1}(h_1) \sigma (h_2 \otimes g_2)
u(h_3g_3)$. 
Then an isomorphism can be carried out by a map $a \otimes h 
\mapsto au(h_1) \otimes h_2$.
\end{prop}

For any commutative ${\bf k}$-algebra $R$ we denote ${\frak Gal}_H(R)$ 
the set of isomorphism
classes of  twisted products $R_{\sigma}[H]$. The following
theorem is somewhat surprising. The informal intuition behind it is 
that the action of $H$ is trivial
and there is no problem how to extend it.
\begin{theorem}
${\frak Gal}_H$ is a covariant functor from the category of
commutative ${\bf k}$-algebras to the category of sets.
\end{theorem}

{\bf Proof.} Let $f:R \rightarrow S$ be a morphism of rings.
It is clear that a cocycle $\sigma$ with values in $R$ gives rise
to a cocycle $f \circ \sigma$ with values in $S$.
We define ${\frak Gal}_H (f) (R_{\sigma}[H])= S_{f \circ \sigma}
[H]$. It is well-defined because equivalent cocycles
with values in $R$ have equivalent continuations to $S$.
Indeed, if $u: H \longrightarrow R$ provides an isomorphism
between $R_{\sigma}[H]$ and $R_{\tau}[H]$ as in Proposition
13 then $S_{f \circ \sigma}[H]$ and $S_{f \circ \tau}[H]$
are isomorphic through $f \circ u$. 
Finally, it is apparent that  ${\frak Gal}_H (f \circ g)
= {\frak Gal}_H (f) \circ {\frak Gal}_H (g)$ for any composition
of algebra maps. $\Box$

We might go on treating ${\frak Gal}_H$ as a presheaf in etale
or Zariski topology to construct the associated sheaf.
It would entitle us to consider $H$-Galois extensions 
(non-commutative torsors) of schemes.
But we prefer to stop here for 
we do not have enough interesting examples
to think about.

\subsection{The structure of the vector bundle of algebras.}
  
We study a vector bundle of finite dimensional algebras on
Spec$\; O$ defined by a Hopf-Galois extension.
Our primary goal is to determine which additional structures
the germs and the fibers inherit.

\begin{theorem} 
The structure of an $H$-comodule algebra can be extended to
$U_\eta , U_{(\eta)} , U_{[\eta]}$ (see~\ref{3}) for each $\eta$. 
The following  inclusions
hold: $U_{[\eta]}^{H} \supseteq O/ \eta$, $U_{\eta}^H \supseteq K_{\eta}$, 
$U_{(\eta)}^H \supseteq O_{(\eta)}$.
If $U$ is prime then 
$U_{(\eta)} \supseteq O_{(\eta)}$ is a 
Hopf-Galois extension with central
invariants under this action.
\end{theorem}
{\bf Proof.} 
Since $H$ is finite dimensional a coaction of $H$ is equivalent
to an action of $H^\ast$. Let $S$ be a complement in $O$ 
of a prime ideal $\eta$. $S$ is a central set of invariants in $U$.
The algebra $U_{(\eta)}$ is isomorphic to  a generalized
algebra of quotients $US^{-1}$. Thus, for $h \in H^\ast$ one
can define
$h \cdot (xs^{-1}) = (h \cdot x)s^{-1}$
 It is clear that $OS^{-1} \cong O_{(\eta)} \subseteq U_{(\eta)}^H$.
On the other hand, given $y \in U_{(\eta)}^H$ one can find
$x \in U, s \in S$ such that $y=xs^{-1}$. Thus, 
$\varepsilon (h) xs^{-1} = (h \cdot x)s^{-1}$ which implies
$\varepsilon (h) x = h \cdot x$ if $S$ is a set of non-zerodivizors of $U$ 
which holds for prime $U$. It implies that $x \in O$ and, therefore,
$y \in O_{(\eta)}$.

Having proved $U_{(\eta)}^H=O_{(\eta)}$, the extension $U_{(\eta)}
 \supseteq O_{(\eta)}$ is $H$-Galois because 
each element of $U_{(\eta )} \otimes H$ 
can be written with a common denominator as $\sum_i a_is^{-1} \otimes h_i$.
We can find
$\sum_i b_i \otimes_O c_i \in U \otimes_O U$ such that
can$( \sum b_i \otimes_O c_i) = \sum a_i \otimes h_i$. It is clear
that  can$( \sum b_is^{-1} \otimes_O c_i) = \sum a_is^{-1} \otimes h_i$.

The other two algebras can be realized as quotients
of $H$-comodule algebras by ideals generated by invariants:
$U_{\eta} \cong U_{(\eta )} / U_{(\eta)}m_\eta$ and
$U_{[\eta]} \cong U / U \eta$. Thus, they admit a canonical
extension of $H^{\ast}$-action. $\Box$ 

\begin{theorem}
If the extension $U \supseteq O$ is cleft then the extensions
$U_{(\eta)} \supseteq O_{(\eta)}$ ,
$U_{[\eta]} \supseteq O/ \eta$, and $U_{\eta} \supseteq K_{\eta}$
are cleft Hopf-Galois with central invariants.
\end{theorem}

{\bf Proof.}
Let $\gamma :H \longrightarrow U$ be a splitting map.
The composition  $H \stackrel{\gamma}{\longrightarrow} O \longrightarrow
O_{(\eta)}$
splits  $U_{(\eta )}$.
Thus, $U_{(\eta )} \cong  U_{(\eta )}^H \otimes H$ as  $O_{(\eta )}$-modules.
On the other hand, changing scalars in $U \cong O \otimes H$ gives
$U_{(\eta )} \cong  O_{(\eta )} \otimes H$. We emphasize that the subspace $H$
of $U$ is the same in the both cases giving a base of $U$ over two subalgebras
$U_{(\eta )}^H \supseteq   O_{(\eta )}$. 
This proves that $U_{(\eta )}^H  =  O_{(\eta )}$
and $U_{(\eta )} \supseteq O_{(\eta)}$ is a cleft $H$-Galois
extension.
The proofs for $U_{\eta}$ and $U_{[\eta]}$ are similar. $\Box$

Another important feature of Hopf-Galois extensions is a Frobenius form \cite{krt}.
Let $\Lambda$ be a left integral
of $H^{\ast}$.
Given an $H$-Galois extension $A \supseteq B$,
we can construct a $B$-bilinear
form $s :A \times A \longrightarrow B$ as
$s(x,y)= x_0y_0 \Lambda (x_1y_1)$ for each $x,y \in A$.
The form $s$ is non-degenerate according to \cite{krt}.

\begin{theorem}
The algebra $U_{\chi}$ is Frobenius for each 
$\chi \in \mbox{Spec}\, O$.
$U_\chi$ is symmetric
if $H$ is unimodular with the antipode
of order 2.
\end{theorem}
{\bf Proof.}
The form $s:U \times U \longrightarrow O$ is non-degenerate \cite{krt}.
By the proof of theorem 6 all specializations $s_{\eta}$ are
non-degenerate proving the first statement.
It was proved in \cite{obe} that $H$ is unimodular
with the antipode of order 2 if and only if
$\Lambda (xy) = \Lambda (yx)$ for each $x,y \in H$.
This implies the second statement. $\Box$

\subsection{Irreducible representations.}

Given a simple  $U$-module $M$,
we get a point 
$\Omega_O (M) \in \mbox{Max} \: O$.
Clearly $M$ can be treated as an $U_{\Omega_O (M)}$-module.
$U_{\Omega_O (M)}$ being a finite-dimensional Frobenius
algebra prompts that the study of blocks and projective covers
of simple modules may be interesting. 

\begin{dfn} Two irreducible $U$-modules $M$ and $N$ are said
to belong to the same block if $\Omega_O (M) = \Omega_O (N)$
and $M$ and $N$ belong to the same block as $U_{\Omega_O (M)}$-modules.
\end{dfn}

It is well-known that the action 
of $Z(U_{\chi})$, the center
of $U_{\chi}$, distinguishes the blocks. On the other hand, rad$\: Z(U_{\chi})$
lies in the kernel of an irreducible representation. These facts
amount to the following lemma.

\begin{lemma}
1. If $M$ and $N$ belong to the same block then $\Omega_Z (M) = \Omega_Z (N)$.

2.  The following two conditions are equivalent for $\chi \in \mbox{\em Max} \: O$:

1) $N$ and $M$ belong to the same block if and only if $\Omega_Z (M) = \Omega_Z (N)$.

2) The composition map $Z \rightarrow Z(U_{\chi}) \rightarrow Z(U_{\chi})/ 
\mbox{\em rad} \: Z(U_{\chi})$ is onto. 
\end{lemma}

If $H$ is copointed (i.e. $H^{\ast}$ is pointed) then little more
is known about simple $U_{\chi}$-modules. The following
fact has been proved in \cite{sch,sc2}. The second statement follows
from the observation that $U_{\chi}^H=K_{\chi}$ has a single simple
module in the cleft case.
 
\begin{prop}
Let us denote the number of simple $A$-modules $s(A)$.
If $H$ is copointed then  $s(U_{\chi}) \leq s(U_{\chi}^{H})s(H)$.
If, furthermore, $U \subseteq O$ is cleft then $s(U_{\chi})$
divides $s(H)$.
\end{prop}

We now assume that the extension $U \supseteq O$ is cleft.  
Let $M$ and $N$ be  two $U_\chi$ modules. According to \cite{lor}
the spectral sequence $E^{p,q}_2= \mbox{Ext}_H^p({\bf k}, \mbox{Ext}_{K_\chi}^q
(M,N))$ converges to $\mbox{Ext}_{U_\chi}^n(M,N)$. Since $K_\chi$ is a field
the elements of the spectral sequence are trivial unless $q=0$.
We state this as a proposition.

\begin{prop}
If $U \supseteq O$ is cleft then for
each $U_{\chi}$-modules $M$ and $N$ we have an isomorphism 
{\em Ext}$^n_{U_{\chi}} (M,N) \cong
\mbox{\em Ext}^n_H({\bf k}, \mbox{\em Hom}_{K_{\chi}}(M,N))$
for each $n$.
\end{prop}


If there is any reason for the blocks of $U_{\chi}$-modules to behave
uniformly with respect to $\chi$ we should be able to find it
looking at the centers of $U_{\chi}$. 
We can succeed showing that the dimensions of the centers
of different $U_{\chi}$ are equal under some strong restrictions
on $U$. 
We assume that $H$ is cocommutative and $U \supseteq O$ is cleft
for the rest of the section. We call a splitting map 
$\gamma :H \longrightarrow U$ equivariant if 
$\gamma (h_1 g S(h_2)) = \gamma (h_1) \gamma (g) \gamma^{-1} (h_2)$
for each $h,g \in H$ .

\begin{lemma} The splitting map is equivariant if and only if
the corresponding cocycle $\sigma$ possesses the following property
$\sigma (a \otimes b) =  \sigma (a_1 b S(a_2) \otimes a_3)$
for each $a,b \in H$.
\end{lemma}
{\bf Proof.}
 It is known that $\sigma (h \otimes g) = 
\gamma (h_1) \gamma (g_1) \gamma^{-1} (h_2g_2)$.
If $\gamma$ is equivariant then
$$\sigma (a_1b_1S(a_2) \otimes a_3)  =
\gamma (a_1b_1S(a_4)) \gamma (a_5) \gamma^{-1} (a_2b_2S(a_3)a_6)=$$
$$\gamma (a_1) \gamma (b_1) \gamma^{-1} (a_2) \gamma (a_3) \gamma^{-1} (a_4b_2S(a_5)a_6)=
\gamma (a_1) \gamma (b_1) \gamma^{-1} (a_2b_2) = \sigma (a \otimes b)$$
If $\sigma$ possesses the property then
$$\gamma (a_1bS(a_2)) = \gamma (a_1bS(a_2)) \gamma (a_3) \gamma^{-1} (a_4) = 
\sigma (a_1b_1S(a_4) \otimes a_5) \gamma (a_2b_2S(a_3)a_6) \gamma^{-1} (a_7) =$$ 
$$\sigma (a_1b_1S(a_2) \otimes a_3) \gamma (a_4b_2S(a_5)a_6) \gamma^{-1} (a_7) =
\sigma (a_1 \otimes b_1) \gamma (a_2b_2) \gamma^{-1} (a_3) =
\gamma (a_1) \gamma (b) \gamma^{-1} (a_2)$$
 $\Box$

An algebra splitting is a trivial example of an equivariant
splitting. It means that the cocycle is $\sigma (x \otimes y)=
\varepsilon (x) \varepsilon (y)1$ and the twisted product
is the tensor product in this case. We provide an example
of a non-trivial equivariant splitting in Section 3.1.
\begin{dfn} We call a linear map $\alpha : H \otimes H \longrightarrow U$
equivariant if for each $a,b \in H$ the following identity holds
$$\alpha (a \otimes b) = \alpha (a_1 b S(a_2) \otimes a_3)$$
\end{dfn}

\begin{lemma} The inverse map under the convolution 
of an equivariant map is equivariant.
The convolution of equivariant maps is equivariant.
\end{lemma}
{\bf Proof.} Define $\beta (a \otimes b) = \alpha^{-1}(a_1 b S(a_2) \otimes a_3)$.
Since $H$ is cocommutative we yield $\beta \ast \alpha (a \otimes b) =
 \beta (a_1 \otimes b_1) \alpha (a_2 \otimes b_2)=
\alpha^{-1}(a_1 b_1 S(a_2) \otimes a_3) \alpha (a_4 b_2 S(a_5) \otimes a_6)=$
$$\alpha^{-1}(a_1 b_1 S(a_4) \otimes a_5) \alpha (a_2 b_2 S(a_3) \otimes a_6)=
\varepsilon (a_1 b_1 S(a_2)) \varepsilon (a_3) = \varepsilon (a) \varepsilon (b)$$
Thus, $\beta = \alpha^{-1}$ proving the first statement.
If $\alpha$ and $\beta$ are equivariant then
$$\alpha \ast \beta (a_1 b S(a_2) \otimes a_3)=
\alpha (a_1 b_1 S(a_4) \otimes a_5) \beta (a_2 b_2 S(a_3) \otimes a_6)=$$
$$\alpha (a_1 b_1 S(a_2) \otimes a_3)\alpha (a_4 b_2 S(a_5) \otimes a_6)=
\alpha (a_1 \otimes b_1) \beta (a_2 \otimes b_2) = \alpha \ast \beta
(a \otimes b) \ \ \ \ \Box$$

\begin{lemma} 
Let $\tau , \pi : H \times H \longrightarrow A$ be two cocycles such that
$\tau \ast \pi^{-1}$ be equivariant. 
Then $a_i \otimes h^i= \sum_i a_i \otimes h^i$ is central in $A_{\tau}[H]$ if and
only if $a_i \otimes h^i$ is central in $A_{\pi}[H]$.
\end{lemma}

{\bf Proof.} Let $\alpha = \tau \ast \pi^{-1}$.
$a_i \otimes h^i$ being central in $A_{\tau}[H]$ provides
$$ a_ib \tau ((h^i)_1 \otimes g_1) \otimes (h^i)_2 g_2 \otimes g_3=
ba_i \tau (g_1 \otimes (h^i)_1)\otimes g_2 (h^i)_2 \otimes g_3$$
for each $b \in A$, $g \in H$. Apply $\rho \otimes \Delta$
to this equality. 
$$a_ib \tau ((h^i)_1 \otimes g_1) \otimes (h^i)_2 g_2 \otimes (h^i)_3 g_3 
\otimes g_4 \otimes g_5=
ba_i \tau (g_1 \otimes (h^i)_1)\otimes g_2 (h^i)_2 \otimes
g_3 (h^i)_3 \otimes g_4 \otimes g_5$$
Apply $\mbox{Id} \otimes \mbox{Id} \otimes 
(m_H \circ (\mbox{Id} \otimes S)) \otimes \mbox{Id}$ where $m_H:H \otimes H
\longrightarrow H$ is the multiplication.
$$a_ib \tau ((h^i)_1 \otimes g_1) \otimes (h^i)_2 g_2 \otimes (h^i)_3 \otimes g_3 
=
ba_i \tau (g_1 \otimes (h^i)_1)\otimes g_2 (h^i)_2 \otimes
g_3 (h_i)_3 S(g_4) \otimes g_5$$
Rewrite using cocommutativity of $H$.
$$(h^i)_1 \otimes g_1 \otimes a_ib \tau ((h^i)_2 \otimes
g_2) \otimes (h^i)_3 g_3=
g_1 (h^i)_1 S(g_2) \otimes g_3 \otimes ba_i 
\tau (g_4 \otimes (h^i)_2)\otimes g_5 (h^i)_3$$
Apply $(m_A \circ (\alpha \otimes \mbox{Id})) \otimes \mbox{Id}$.
$$\alpha((h^i)_1 \otimes g_1) a_ib \tau ((h^i)_2 \otimes g_2) 
\otimes (h^i)_3 g_3=
\alpha (g_1 (h^i)_1 S(g_2) \otimes g_3) ba_i 
\tau (g_4 \otimes (h^i)_2)\otimes g_5 (h^i)_3$$
Since $A$ is commutative and $\alpha$ is equivariant it may be rewritten as
$$a_ib \pi ((h^i)_1 \otimes g_1) \otimes (h^i)_2 g_2=
ba_i \pi (g_1 \otimes (h^i)_1)\otimes g_2 (h^i)_2$$
which proves that $a_i \otimes h^i$ is central in $A_{\pi}[H]$. $\Box$

\begin{theorem} 
Let us assume that the extension $U \supseteq O$ admits equivariant
splitting and $H$ is cocommutative.
Then the dimension of $Z(U_{\chi})$, the center of $U_{\chi}$,
over $K_{\chi}$ does not
depend on $\chi$.
\end{theorem}

{\bf Proof.}
Any two fields $K_\chi$ and $K_\eta$ have a common overfield $F$.
Let $e_i$ be a basis of $F$ over $K_\chi$. An element of the center
of $U_{\chi} \otimes_{K_{\chi}} F$ may be written as $\sum_i a_i \otimes
e_i$ with $a_i \in U_{\chi}$. Since it commutes with elements
of $U_{\chi} \otimes 1$ the equality $Z(U_{\chi} \otimes_{K_{\chi}}F)=
Z(U_{\chi}) \otimes_{K_{\chi}}F$ holds and the dimension of the center
does not change with the field extension. Both $U_{\chi} \otimes_{K_{\chi}}F 
\supseteq F$ and $U_{\eta} \otimes_{K_{\chi}}F \supseteq F$ inherit
equivariant splittings from $U \supseteq O$ and the theorem follows
from the lemmas above. $\Box$

\section{Examples.}

 A number of interesting examples of central Hopf-Galois
extensions arise in study of Hopf algebras. 
Given a Hopf algebra $U$ with central Hopf subalgebra $O$,
we  treat $U \supseteq O$ as an $H$-Galois extension
with $H=U/U(O \cap \mbox{ker} \; \varepsilon )$. 
The map $\rho :U \longrightarrow U \otimes H$
is $\rho (x) =x_1 \otimes \bar{x}_2$.
$H$ is finite dimensional if and only if the module
$U_O$ is finitely generated. 


We discuss one more feature before doing examples. It is 
the winding action. It describes all $U_{\chi}$ isomorphic
to $H$. This has been shown in \cite[Exercise 5.3.4]{far}
for the reduced enveloping algebras (see Section 3.2).

\begin{prop} $U_{\chi}$ is isomorphic to $H$ if and only if
$U_{\chi}$ has a one-dimensional representation.
\end{prop}

{\bf Proof.}
The direct implication is apparent. Let $U_{\chi}$ have 
a one-dimensional representation $\alpha :U_{\chi} \longrightarrow
{\bf k}$ for some $\chi$. This particularly means
${\bf k}=K_{\chi}$. The composition 
$\tilde{\alpha} :U \longrightarrow U_{\chi} 
\stackrel{\alpha}{\longrightarrow} {\bf k}$ is a one-dimensional
representation of $U$.
We define $\hat{\alpha} :U \longrightarrow U$ as $\hat{\alpha} (x)=
\tilde{\alpha} (x_1)x_2$. Clearly, $U/I \cong U/ \hat{\alpha}(I)$
for each ideal $I$. $\tilde{\alpha} (x) = \tilde{\alpha} (x_1) \varepsilon (x_2) =
\varepsilon ( \hat{\alpha} (x))$ proving that
$\hat{\alpha}( \mbox{ker} \; \tilde{\alpha})= \mbox{ker} \; \varepsilon$.
Finally, $\hat{\alpha} (O)=O$ since $O$ is a Hopf subalgebra and,
therefore, $\hat{\alpha} (U(O \cap \mbox{ker} \; \tilde{\alpha})) =
U(O \cap \mbox{ker} \; \varepsilon )$. $\Box$

\subsection{Groups.}

Let $H$ be a central subgroup of a group $G$.
Then any set splitting $\tilde{\gamma} :G/H \longrightarrow G$
gives rise to a coalgebra splitting $\gamma :{\bf k}G/H
\longrightarrow {\bf k}G$. This example fits precisely
to the setup of the present paper if $H$ is finitely
generated of finite index.
Let us further assume that $G$ is a free Abelian group
and $H$ is proper. Then there exists no Hopf algebra
splitting $\gamma :{\bf k}G/H \longrightarrow {\bf k}G$
for $G/H$ must be a subgroup of $G$ otherwise.
This provides an example of a non-trivial equivariant
splitting because 
$\gamma (h_1 g S(h_2)) = \gamma (gh_1 S(h_2)) = \gamma (g)
\varepsilon (h) = 
\gamma (g) \gamma (h_1) 
\gamma^{-1} (h_2) = \gamma (h_1) \gamma (g) \gamma^{-1} (h_2)$.

\subsection{Restricted Lie algebras.}

 Let $L$ be a finite dimensional restricted Lie algebra
over an algebraically closed field {\bf k} of characteristic $p$.
The universal enveloping algebra $U(L)$ has a central Hopf
subalgebra $O$ generated by $x^p-x^{[p]}$ for $x \in L$.
The quotient $U(L)/U(L) ( \mbox{ker} \varepsilon \cap O)$ is the
restricted enveloping algebra $u(L)$.

For any $\chi \in L^{\ast}$ the reduced enveloping algebra $U^{\chi}$
is a quotient of $U(L)$ by the ideal generated by $x^p-x^{[p]}- 
\chi (x)^p1_{U(L)}$ for all $x \in L$ \cite{far,fri}.  
We should notice that the standard notation $U_{\chi}$ for
the reduced enveloping algebra has already been used
for a different object 
in the present paper. We will see pretty soon that our algebras $U_{\chi}$
are not far from the reduced enveloping algebras.
The map $x \mapsto x^p-x^{[p]}$ can be extended to
a semilinear isomorphism between the symmetric algebra
$S(L)$ and $O$. This defines an isomorphism
between Spec$\: O$ and the Frobenius shift $L^{\ast \; (1)}$. 
The fiber algebras at the closed points are 
precisely the reduced enveloping algebras we have just defined.

\begin{lemma}
$U^{\chi}$ is isomorphic to $U_{\chi^p}$. If we think of $U^{\chi}$
as a bundle of algebras over $L^{\ast}$ then the isomorphism
is a pull-back along the map $L^{\ast} \stackrel{F}{\longrightarrow}
L^{\ast \; (1)} \stackrel{\cong}{\longrightarrow} \mbox{Spec} \; O$.
\end{lemma}

Furthermore, $U(L) \supseteq O$ is a cleft $u(L)$-Galois extension
and, therefore,
all algebras $U_{\chi}$ are twisted products. Its being cleft
follows from $u(L)$ being pointed
but we can also construct the splitting explicitly.

\begin{prop}
Let us choose an ordered basis $e_i$ of $L$. It gives rise
to PBW-bases on $U(L)$ and $u(L)$. Denoting their elements
as $e^\alpha$ we can construct a splitting map
$\gamma : u(L) \longrightarrow U(L)$ by $\gamma (e^{\alpha})
=e^{\alpha}$.
\end{prop}

{\bf Proof.} It is apparent that $\gamma$ is a linear splitting
as well as a coalgebra map. $\Box$

As an application of this fact we may give an explicit formula
for multiplication in $U_{\chi}$. Having chosen an ordered basis of $L$,
we identify $U_{\chi}$ and $u(L)$ as vector spaces 
using a PBW-basis in both. We use the standard notation
for the Hopf algebra structure maps of $u(L)$;
 the multiplication of $U_\chi$ is denoted $\circ$.
We think of $\chi$ as a linear map from $O$ to ${\bf k}$.
The map $\gamma : u(L) \longrightarrow U(L)$ is constructed in Proposition 29.
The following proposition easily follows from 
the discussion above as well as the observation that $\gamma^{-1}(x)=
S( \gamma (x))$.

\begin{prop}
The element $\sigma (x \otimes y)= \gamma (x_1) \gamma (y_1)
S( \gamma (x_2y_2))$ belongs to $O$ for each $x,y \in u(L)$.
The multiplication of $U_\chi$ can be written as
$x \circ y = \chi ( \sigma (x_1 \otimes y_1)) x_2y_2$.
\end{prop}

The case of $L=sl_2({\bf k})$ is the easiest to visualize.
We assume that $p > 2$.
Let $e,f,h$ be the standard basis of $sl_2({\bf k})$. Following
\cite{rud}, we denote $x=e^p,
y=f^p, z=h^p-h, t=(h+1)^2-4ef$. Then $O$ is isomorphic
to the polynomial algebra ${\bf k}[x,y,z]$ and $Z$ is isomorphic
to ${\bf k}[x,y,z,t]/ \langle t^p-2t^{(p+1)/2} +t-(z^2-4xy) \rangle$.

We  treat points of Max$\; O$ as elements of $L^{\ast \; (1)}$
as well as triples $(x,y,z)$. Points of Max$\; Z$ are quadruples
$(x,y,z,t)$ satisfying the equation: 
\begin{equation} t^p-2t^{(p+1)/2} +t=(z^2-4xy) \label{4} \end{equation}
One may observe three types of points in $L^{\ast (1)}$.
If $\chi =(x,y,z)$ satisfies $z^2-4xy \neq 0$ then the equation~\ref{4} on $t$
has $p$ distinct roots and $U_{\chi}$ is a direct sum of $p$ copies of $M_p({\bf k})$. Over the cone
$z^2-4xy = 0$ the equation~\ref{4} has 1 single root 0 and $\frac{p-1}{2}$
double roots $i^2 \; $mod$\; p$, $i= 1, \ldots , \frac{p-1}{2}$. 
A point $(x,y,z,t)$ is singular if and only if $x=y=z=0$ and $t \neq 0$
as easy to see by a direct computation.
The point $(0,0,0,0)$ corresponds to the Steinberg module for $u(sl_2({\bf k}))$.
A singular point $(0,0,0,i)$ corresponds to the block consisting
of 2 restricted $sl_2({\bf k})$-modules of dimensions $i$ and $p-i$.
For remaining $\chi$ of the cone
all $U_{\chi}$ are isomorphic. This algebra has $\frac{p+1}{2}$
simple modules of dimension $p$ \cite{fri}. One of them is projective;
the others comprise distinct blocks and have double projective covers.

One may consult \cite{fri,fr2,kac} for the case of any
classical semisimple Lie algebra.  However, the precise details
are unknown in the general case.
Another case which may be treated explicitly is one of a completely
solvable Lie algebra $L$ \cite{far,kac}.

\subsection{Quantum groups.}
 The quantum enveloping algebra $U_q (g)$ of a semisimple Lie algebra $g$
is generated by the
elements $E_\alpha , F_\alpha , K_i$. If $q$ is a primitive $l$-th root
of unity then $O$ is the subalgebra generated by $E_{\alpha}^l ,
 F_{\alpha}^l , K_i^l$. If $n$ is the dimension of $g$ and
$r$ is the rank
we get a bundle of algebras on
${\Bbb C}^{n-r} \times 
{\Bbb C}^{r}$. This resembles the case of restricted Lie algebras.
The construction
has been developed  in \cite{dck}.

\subsection{Quantum linear groups.}

 The following example has been worked out for any
complex semisimple group in \cite{dcl}.
We follow \cite{par} with an elementary approach to $SL (n)$.
The quantum linear group $SL_q (n)$ is a non-commutative $\mathbb C$-algebra
generated by $X_{ij}, \ i,j = 1, \ldots , n$. The following are defining
relations:
$$X_{ri}X_{rj}=q^{-1}X_{rj}X_{ri} \mbox{ if } i < j$$
$$X_{ri}X_{si}=q^{-1}X_{si}X_{ri} \mbox{ if } r<s$$
$$X_{ri}X_{sj}=X_{sj}X_{ri} \mbox{ if } i > j \mbox{ and } r<s$$
$$X_{ri}X_{sj}-X_{sj}X_{ri}=(q^{-1}-q)X_{si}X_{rj} \mbox{ if } i < j 
\mbox{ and } r<s$$
$$\sum_{\sigma \in S_n} (-q)^{-l( \sigma )}X_{1 \sigma (1)} \cdots
X_{n \sigma (n)} =1$$
 $l ( \sigma )$ denotes the length of a permutation $\sigma$ which is
the minimal number of transpositions required to represent $\sigma$.
$SL_q (n)$ has a structure of a Hopf algebra with the comultiplication
$\Delta (X_{ij}) = \sum_k X_{ik} \otimes X_{kj}$ and the counity
$\varepsilon (X_{ij})= \delta_{ij}$. The definition of the antipode 
requires algebraic complements; the inquisitive reader may refer to \cite{par}.

If $q$ is an $m$-th root of unity the subalgebra $O$ generated by $X_{ij}^m$
is central. Furthermore, it is a Hopf subalgebra isomorphic to ${\mathbb C}[SL_n]$.
Thus, following the construction of the present paper one gets a bundle of algebras
on $SL_n (\mathbb C)$. 
If $q$ is a primitive root and $m$ is odd then 
then there are at most $(n!)^2$ non-isomorphic fiber algebras
in this bundle. Subsets over which the fiber algebras are isomorphic
are parametrized by pairs $(w_1,w_2) \in S_n \times S_n$ in \cite{dcl}.
We think of $S_n$ as a subset of $SL_n(\mathbb C)$ using a splitting
of the map $N \longrightarrow S_n$ where $N$ is the normalizer
of the subgroup of diagonal matrices. 
 If $B_+$ and $B_-$ are the subgroups of upper and lower
triangular matrices then the subset corresponding to $(w_1,w_2)$
is $(B_+w_1B_+) \cap (B_-w_2B_-)$.

\subsection{Subvarieties}
Restricting the bundle of algebras to a locally closed subset of Spec$\; O$
provides an example of a Hopf-Galois extension of not necessarily
affine spaces.
For instance, a restriction to nilpotent orbits of a modular semisimple 
Lie algebra  of the bundle in Example 3.2 has been studied in \cite{fr2}.

\subsection{Quantum commutative algebras.}

The following idea has been utilized in \cite{coh} to get examples
of central Hopf-Galois extensions. Let $(H, \langle , \rangle )$
be a finite-dimensional coquasitriangular Hopf algebra.
We consider an $H$-Galois extension $U \supseteq O$ such that
$U$ is quantum commutative; i.e. 
$xy = y_0 x_0 \langle x_1 , y_1 \rangle$ for each $x,y \in U$. 
These conditions imply that the subalgebra $O$ is central.
Indeed, $\langle 1,h \rangle = \varepsilon (h)$ for 
every $h \in H$ \cite[10.1.8]{mon}.
Given $x \in U, y \in U$, we have $xy= y_0 x \langle 1 , y_1 \rangle
=yx$.

\begin{center}
{\bf Address}

Department of Mathematics \& Statistics, 

LGRT, UMass, Amherst, MA, 01003.

{\bf E-mail:} rumynin@@math.umass.edu
\end{center}
\end{document}